\documentclass[amssymb, 11pt]{amsart}
\usepackage{latexsym}

\newlength{\standardunitlength}
\newtheorem{prop}{Proposition}[section]

\newtheorem{lemma}[prop]{Lemma}
\newtheorem{cor}[prop]{Corollary}
\newtheorem{theorem}[prop]{Theorem}

\begin{document}

%\begin{center}
\title [Card shuffling] {Card shuffling and the decomposition of tensor products}
%\end{center}

\author{Jason Fulman}
\address{University of Pittsburgh\\
Pittsburgh, PA}
\email{fulman@math.pitt.edu}

\keywords{Card shuffling, tensor product, Markov chain, Robinson-Schensted-Knuth correspondence}

\subjclass{}

\date{July 1, 2003}

\thanks{Fulman was partially supported by National Security Agency
grant MDA904-03-1-0049.}

\begin{abstract} Let $H$ be a subgroup of a finite group $G$. We use Markov chains to quantify how large $r$ should be so that the decomposition of the $r$ tensor power of the representation of $G$ on cosets on $H$ behaves (after renormalization) like the regular representation of $G$. For the case where $G$ is a symmetric group and $H$ a parabolic subgroup, we find that this question is precisely equivalent to the question of how large $r$ should be so that $r$ iterations of a shuffling method randomize the Robinson-Schensted-Knuth shape of a permutation. This equivalence is rather remarkable, if only because the representation theory problem is related to a reversible Markov chain on the set of representations of the symmetric group, whereas the card shuffling problem is related to a nonreversible Markov chain on the symmetric group. The equivalence is also useful, and results on card shuffling can be applied to yield sharp results about the decomposition of tensor powers. \end{abstract}

\maketitle

\section{Introduction} \label{intro}

	Let $\chi$ be a faithful character of a finite group $G$. A well known theorem of Burnside and Brauer (\cite{I}) states that if $\chi(g)$ takes on exactly $m$ distinct values for $g \in G$, then every irreducible character of $G$ is a constituent of one of the characters $\chi^j$ for $0 \leq j < m$. It is very natural to investigate the decomposition of $\chi^j$, and the results in this paper are a step in that direction.

	Let $Irr(G)$ denote the set of irreducible representations of a finite group $G$. The Plancherel measure on $Irr(G)$ is a probability measure which assigns mass $\frac{dim(\rho)^2}{|G|}$ to $\rho$. The symbol $\chi^{\rho}$ denotes the character associated to the representation $\rho$. The notation $Ind,Res$ stands for induction and restriction of class functions. We remind the reader that the character of the $r$-fold tensor product of a representation of $G$ is given by raising the character to the rth power. The inner product $<f_1,f_2>$ denotes the usual inner product on class functions of $G$ defined by \[ \frac{1}{|G|} \sum_{g \in G} \overline{f_1(g)} f_2(g).\] Thus if $f_1$ is an irreducible character and $f_2$ any character, their inner product gives the multiplicity of $f_1$ in $f_2$. We let $g^G$ denote the conjugacy class of $g$ in $G$.

	In Section \ref{general} of this paper, we prove the following result.

\begin{theorem} \label{main1} Let $H$ be a subgroup of a finite group $G$. Let $\pi$ denote the Plancherel measure of $G$. Suppose that $|G|>1$. Let \[ \beta = max_{g \neq 1} \frac{|g^G \cap H|}{|g^G|} = \frac{|H|}{|G|} max_{g \neq 1} Ind_H^G(1)[g].\]  Then \[   \sum_{\rho \in Irr(G)} |(\frac{|H|}{|G|})^r dim(\rho) <\chi^{\rho}, (Ind_{H}^G(1))^r > - \pi(\rho)| \leq |G|^{1/2} \beta^r.\] \end{theorem} Note that if $\beta<1$, the right hand side approaches 0 as $r \rightarrow \infty$. The quantity $\beta$ has been carefully studied in the (most interesting) case that $G$ is simple and $H$ a maximal subgroup; references and an example where $H$ is not maximal are given in Section \ref{general}.

	The idea behind the proof of Theorem \ref{main1} is to investigate a natural Markov chain $J$ on the set of irreducible representations of $G$. This chain is essentially a probabilistic reformulation of Frobenius reciprocity. This chain can be explicitly diagonalized and then Theorem \ref{main1} follows from spectral theory of reversible Markov chains. In fact Theorem \ref{main1} is a generalization of a result in our earlier paper \cite{F1}, where this Markov chain arose for the symmetric group case $H=S_{n-1}$ and $G=S_n$ and was combined with Stein's method to sharpen a result of Kerov on the asymptotic normality of random character ratios of the symmetric group on transpositions.

	 The main insight of the current paper is that when $G$ is the
	 symmetric group $S_n$ and $H$ is a parabolic subgroup, the
	 bound of Theorem \ref{main1} can be improved by card
	 shuffling. Let us describe this in detail for the case
	 $H=S_{n-1}$. In Theorem \ref{main1}, $\beta=1-\frac{2}{n}$,
	 and one can see using Stirling's approximation for $n!$ that for
	 $r>\frac{nlog(n)+c}{2log(\frac{1}{\beta})}$, the bound in
	 Theorem 1 is at most $(2 \pi)^{1/4} e^{-c}$ (and hence
	 small). Thus $r$ slightly more than $\frac{n^2 log(n)}{4}$
	 suffices to make the bound small. The bound of Theorem \ref{main1} is proved by
	 analyzing a certain Markov chain $J$ on $Irr(S_n)$, started
	 at the trivial representation. The irreducible representations
	 of $S_n$ correspond to partitions of $n$ (the one row
	 partition is the trivial representation), so $J$ is a Markov
	 chain on partitions. We show that the distribution on
	 partitions given by taking $r$ steps according to $J$ has a
	 completely different description. Namely starting from the
	 identity permutation (viewed as $n$ cards in order), perform
	 the following procedure $r$ times: remove the top card and
	 insert it into a uniformly chosen random position. This gives
	 a nonuniform random permutation, and there is a natural map
	 called the Robinson-Schensted-Knuth or RSK correspondence
	 (see \cite{Sa} for background), which associates a partition
	 to a permutation. We will show that applying this
	 correspondence to the permutation obtained after $r$
	 iterations of the top to random shuffle gives {\it exactly}
	 the same distribution on partitions as that given by $r$ iterations of the chain $J$ started at the trivial representation. This will allow us to use facts about card
	 shuffling to sharpen to $\frac{n^2log(n)}{4}$ to roughly
	 $nlog(n)$, and even to see that the $nlog(n)$ is 
	 sharp to within a factor of two. Precise statements and results for more general
	 parabolic subgroups are given in Section \ref{symmetric}.

	To close the introduction we make some remarks. First, recall
 that a Markov chain $M$ on a finite set $X$ is called reversible with
 respect to the probability measure $\mu$ on $X$ if $\mu(x) M(x,y) =
 \mu(y) M(y,x)$ for all $x,y$ (it follows that $\mu$ is
 stationary for $M$, i.e. that $\mu(y) = \sum_x \mu(x) M(x,y)$ for all $y$). The top to random shuffle and its cousins which
 arise in connection with parabolic subgroups are nonreversible
 chains. Thus it is rather miraculous that the top to random shuffle
 has real eigenvalues; this observation is the starting point of a
 general theory \cite{BHR}. And it is doubly surprising that the top
 to random shuffle should be connected with the reversible chains
 $J$. Second, the problem of studying the convergence rate of the RSK
 shape after iterated shuffles to the RSK shape of a random
 permutation is of significant interest independent of its application in
 this paper. It is closely connected with random matrix theory and in
 some cases with Toeplitz determinants. See \cite{St}, \cite{F2},
 \cite{F3} and the references therein for details. Third, since
 Solomon's descent algebra generalizes to finite Coxeter groups, it is
 likely that the results in this paper can be pushed through to that
 setting. (However that would require an analog of the RSK
 correspondence for finite Coxeter groups).

\section{General groups} \label{general}

	This section proves Theorem \ref{main1} and gives an example. Throughout this section $X=Irr(G)$ is the set of irreducible representations of a finite group $G$, endowed with Plancherel measure $\pi_G$. We also suppose that we are given a subgroup $H$ of $G$.

	To begin, we use $H$ to construct a Markov chain on $G$ which is reversible with respect to $\pi_G$. For $\rho$ an irreducible representation of $G$ and $\tau$ an irreducible representation of $H$, we let $\kappa(\tau,\rho)$ denote the multiplicity of $\tau$ in $Res_H^G(\rho)$. By Frobenius reciprocity, this is the multiplicity of $\rho$ in $Ind_H^G(\tau)$.

\begin{prop} The Markov chain $J$ on irreducible representations of $G$ which moves from $\rho$ to $\sigma$ with probability \[ \frac{|H|}{|G|} \frac{dim(\sigma)}{dim(\rho)} \sum_{\tau \in Irr(H)} \kappa(\tau,\rho) \kappa(\tau,\sigma) \] is in fact a Markov chain (the transition probabilities sum to 1), and is reversible with respect to the Plancherel measure $\pi_G$.
\end{prop}

\begin{proof} First let us check that the transition probabilities sum to 1. Indeed, \begin{eqnarray*}
& & \sum_{\sigma \in Irr(G)} \frac{|H|}{|G|} \frac{dim(\sigma)}{dim(\rho)} \sum_{\tau \in Irr(H)} \kappa(\tau,\rho) \kappa(\tau,\sigma)\\
& = & \frac{|H|}{|G|} \frac{1}{dim(\rho)} \sum_{\tau \in Irr(H)} \kappa(\tau,\rho) \sum_{\sigma \in Irr(G)} dim(\sigma) \kappa(\tau,\sigma)\\
& = & \frac{1}{dim(\rho)} \sum_{\tau \in Irr(H)} \kappa(\tau,\rho) dim(\tau)\\
& = & 1. \end{eqnarray*} The second equality follows since the dimension of a representation induced from a subgroup is its original dimension multiplied by the index of the subgroup. 

	The reversibility with respect to Plancherel measure is immediate from the definitions. \end{proof}

	Next we quickly review some facts from Markov chain theory. We consider the space of real valued functions $\ell^2(\pi)$ with the norm \[ ||f||_2 = \left( \sum_x |f(x)|^2 \pi(x) \right)^{1/2}.\] If $J(x,y)$ is the transition rule for a Markov chain on $X$, the associated operator (also denoted by $J$) on $\ell^2(\pi)$ is given by $Jf(x) = \sum_y J(x,y) f(y)$. Let $J^r(x,y)=J_x^r(y)$ denote the chance that the Markov chain started at $x$ is at $y$ after $r$ steps.

	If the Markov chain with transition rule $J(x,y)$ is reversible with respect to $\pi$ (i.e. $\pi(x) J(x,y) = \pi(y) J(y,x)$ for all $x,y$), then the operator $J$ is self adjoint with real eigenvalues \[ -1 \leq \beta_{min}=\beta_{|X|-1} \leq \cdots \leq \beta_1 \leq \beta_0=1.\] Let $\psi_i$ ($i=0,\cdots,|X|-1$) be an orthonormal basis of eigenfunctions such that $J\psi_i = \beta_i \psi_i$ and $\psi_0 \equiv 1$. Define $\beta = max \{|\beta_{min}|,\beta_1\}$.

	The total variation distance between two probability measures $Q_1,Q_2$ on a set $X$ is defined as $||Q_1- Q_2||_{TV} = \frac{1}{2} \sum_{x \in X} |Q_1(x)-Q_2(x)|$. It is elementary that $||Q_1-Q_2||_{TV}=max_{A \subseteq X} |Q_1(A)-Q_2(A)|$. Thus when the total variation distance is small, the $Q_1$ and $Q_2$ probabilities of any event $A$ are close.

	The following lemma is well-known; for a proof see \cite{DS}. 

\begin{lemma} \label{classic}
\begin{enumerate}
\item $2 ||J_x^r - \pi||_{TV} \leq ||\frac{J_x^r}{\pi}-1||_2$.
\item $J^r(x,y) = \sum_{i=0}^{|X|-1} \beta_i^r \psi_i(x) \psi_i(y) \pi(y)$.
\item $||\frac{J_x^r}{\pi}-1||_2^2 = \sum_{i=1}^{|X|-1} \beta_i^{2r} |\psi_i(x)|^2 \leq \frac{1-\pi(x)}{\pi(x)} \beta^{2r}$.
\end{enumerate}
\end{lemma}

\begin{prop} \label{spectrum} Let $G$ be a finite group and $H$ a subgroup of $G$. Then the eigenvalues and eigenfunctions of the operator $J$ are indexed by conjugacy classes $C$ of $G$.
\begin{enumerate}
\item The eigenvalue parameterized by $C$ is $\frac{|C \cap H|}{|C|}$. 
\item The orthonormal basis of eigenfunctions $\psi_C$ is defined by $\psi_C(\rho) = \frac{|C|^{\frac{1}{2}} \chi^{\rho}(C)}{dim(\rho)}$.
\end{enumerate}
\end{prop}

\begin{proof} First, note that the transition probability in the definition of $J$ can be rewritten as follows: \begin{eqnarray*} & & \frac{|H| dim(\sigma)}{|G| dim(\rho)} <\chi^{\sigma},Ind_{H}^G Res^G_H (\chi^{\rho})>\\
& = & \frac{|H| dim(\sigma)}{|G| dim(\rho)} \frac{1}{|G|} \sum_{g \in G} \overline{\chi^{\sigma}(g)} \frac{1}{|H|} \sum_{t \in G \atop t^{-1}gt \in H} \chi^{\rho}(t^{-1}gt)\\
& = & \frac{dim(\sigma)}{dim(\rho)} \frac{1}{|G|} \sum_{g \in G} \overline{\chi^{\sigma}(g)} \chi^{\rho}(g) \frac{|g^G \cap H|}{|g^G|}. \end{eqnarray*} The first equality used the well known formula for induced characters \cite{I}. 

	Now to see that $\psi_C$ is an eigenfunction with the asserted eigenvalue, one calculates that
\begin{eqnarray*}
& & \sum_{\sigma \in Irr(G)} \frac{dim(\sigma)}{dim(\rho)} \frac{1}{|G|} \sum_{g \in G} \overline{\chi^{\sigma}(g)} \chi^{\rho}(g) \frac{|g^G \cap H|}{|g^G|}|C|^{\frac{1}{2}} \frac{\chi^{\sigma}(C)}{dim(\sigma)}\\
& = & \frac{|C|^{\frac{1}{2}} }{dim(\rho)} \sum_{g \in G} \frac{|g^G \cap H|}{|g^G|} \frac{\chi^{\rho}(g)}{|G|} \sum_{\sigma \in Irr(G)} \overline{\chi^{\sigma}(g)} \chi^{\sigma}(C)\\
& = & \frac{|C|^{\frac{1}{2}} \chi^{\rho}(C)}{dim(\rho)} \sum_{g \in C} \frac{|g^G \cap H|}{|g^G|} \frac{1}{|g^G|}\\
& = & \frac{|C|^{\frac{1}{2}} \chi^{\rho}(C)}{dim(\rho)} \frac{|C \cap H|}{|C|}. \end{eqnarray*} Note that the second inequality used the orthogonality relations of the characters of $G$. 

	Finally, the fact that $\psi_C$ are orthonormal follows from the orthogonality relations for irreducible characters. It is well known that $\psi_C$ are a basis. \end{proof}

	Next we prove Theorem \ref{main1} from the introduction.

\begin{proof} First note that the equivalence of the definitions of $\beta$ follows from the general formula for induced characters. Now let $1$ denote the trivial representation of $G$. From Proposition \ref{spectrum} and part 2 of Lemma \ref{classic}, \begin{eqnarray*} J_{1}^r(\rho) & = & dim(\rho) \sum_{C} (\frac{|C \cap H|}{|C|})^r \frac{|C| \chi^{\rho}(C)}{|G|}\\
& = & dim(\rho) \frac{1}{|G|} \sum_{g \in G} (\frac{|g^G \cap H|}{|g^G|})^r \chi^{\rho}(C)\\
& = & dim(\rho) (\frac{|H|}{|G|})^r <\chi^{\rho}, (Ind_{H}^G(1))^r>, \end{eqnarray*} where in the third equality we have used the well known formula for induced characters used in the proof of Proposition \ref{spectrum}. The theorem now follows from part 1 of Proposition \ref{spectrum} and parts 1 and 3 of Lemma \ref{classic}. \end{proof}

{\bf Remarks:}

\begin{enumerate}
\item The quantity $\beta$ has been well studied in the case that $G$ is simple and $H$ is a maximal subgroup of $G$. See for instance \cite{GK}, \cite{LS} and the references therein. We defer discussion of the case that $G=S_n$ and $H$ is a parabolic subgroup to Section \ref{symmetric}. 

\item Observe that if $\beta=1$ the upper bound of Theorem \ref{main1} is useless. And it can happen that $\beta=1$. For instance if $H$ is a nontrivial normal subgroup of $G$, there are conjugacy classes of $G$ contained in $H$. On the representation theory side, suppose for simplicity that $H$ is normal of index 2. Then except in trivial cases, the state space of the Markov chain $J$ isn't connected, so the the quantity bounded in Theorem \ref{main1} won't go to 0 as $r \rightarrow \infty$. Indeed, either $Ind_{H}^G Res_{H}^G (\rho)$ is 2 copies of $\rho$ or else the sum of $\rho$ and the conjugate representation of $\rho$ (page 64 of \cite{FH}).

\item Observe that if $\beta=0$, then $|H|=1$ which implies that the decomposition of $Ind_H^G(1)$ is given exactly by Plancherel measure. Then the bound in Theorem \ref{main1} is an equality.
\end{enumerate}

	To conclude this section we compute $\beta$ in the case that $G=GL(n,q)$ and $H=GL(n-1,q)$ (which is not a maximal subgroup). There are clearly more examples in this direction which can be worked out using Wall's formulas for conjugacy class sizes \cite{W}-though as in Proposition \ref{Glcase} below some (minor) effort is required to determine when $\frac{|g^G \cap H|}{|g^G|}$ is largest for nontrivial $g$. However as we have no need for them we stop here.

\begin{prop} \label{Glcase} Suppose that $G=GL(n,q)$ and $H=GL(n-1,q)$, and that $n \geq 2$. Then $\beta=\frac{(1-1/q^{n-1})}{q^2 (1-1/q^n)}$ for $q>2$ and $\beta=\frac{(1-1/q^{n-2})}{q^2 (1-1/q^n)}$ for $q=2$. 
\end{prop}

\begin{proof} The conjugacy classes $C$ of
 $GL(n,q)$ are parameterized by all ways of associating a partition
 $\lambda_{\phi}$ to each monic irreducible polynomial $\phi(z)$ with coefficients in
 $F_q$ such that $|\lambda_{z}|=0$ and $\sum_{\phi} deg(\phi)
 |\lambda_{\phi}|=n$. Here $|\lambda|$ denotes the size of a partition
 $\lambda$ and $deg(\phi)$ denotes the degree of the polynomial
 $\phi$. Moreover the size of the conjugacy class with this data is (\cite{M}, page 181) \[
 \frac{|GL(n,q)|}{\prod_{\phi} \prod_{j \geq 1} q^{deg(\phi)
 (\lambda_{\phi,j}')^2} (1-1/q^{deg(\phi)}) \cdots (1-1/q^{deg(\phi)
 m_j(\lambda_{\phi})})}.\] Here $m_j(\lambda_{\phi})$ is the number of parts
 of $\lambda_{\phi}$ of size $j$, and $\lambda_{\phi,j}'$ is the
 number of parts of $\lambda_{\phi}$ of size at least $j$. In order that
 $\frac{|g^G \cap H|}{|g^G|}$ is nonzero, it is necessary that $g$ has
 its conjugacy data satisfying $m_1(\lambda_{z-1}(g)) \geq 1$. Then $g^G \cap H$ is a single conjugacy class of $H$, with
 conjugacy data the same as for $g$ except that a part of size 1 is
 removed from the partition corresponding to the polynomial
 $z-1$. Thus one sees that \[ \frac{|g^G \cap H|}{|g^G|} =
 \frac{|GL(n-1,q)|}{|GL(n,q)|} 
 (1-1/q^{m_1(\lambda_{z-1}(g))}) q^{2 \lambda_{z-1,1}'(g)-1}.\]

	Thus to find $\beta$, it is necessary study the maximum of the function \[ (1-1/q^{m_1(\lambda)}) q^{2 \lambda_1'} \]
among partitions $\lambda$ of size at most $n$ having at least 1
	part equal to 1, but excluding the partition of size $n$ which consists of $n$ 1's. Here $m_1(\lambda)$ denotes the number of parts of $\lambda$ of size 1, and $\lambda_1'$ denotes the number of parts of $\lambda$. It is straightforward to see that if $|\lambda|<n$, this function is maximized when $|\lambda|=n-1$ and $\lambda$ consists of $n-1$ 1's. For $|\lambda|=n$ it is straightforward that the function is maximized for the partition consisting of 1 part of size 2 and $n-2$ parts of size 1. Comparing these two cases one sees that the maximum occurs for the first case. The first case occurs for $q>2$ but can not occur for $q=2$ (since $z-1$ is the only polynomial of degree 1 with nonzero constant term), and for $q=2$ it is straightforward to see that the second case is the maximum. \end{proof}

\section{Symmetric groups} \label{symmetric}

	This section considers the Markov chain $J$ in the case of the symmetric group and develops connections with card shuffling. We assume throughout that the reader is familiar with the Robinson-Schensted-Knuth (RSK) correspondence. See \cite{Sa} for background on this topic. 	

	Consider the symmetric group $S_n$. Let $\Pi=\{\epsilon_1-\epsilon_2, \cdots, \epsilon_{n-1}-\epsilon_n\}$ be a set of simple roots for the root system consisting of the $n(n-1)$ vectors $\epsilon_i-\epsilon_j$, where $1 \leq i \neq j \leq n$. The positive roots are $\epsilon_i - \epsilon_j$ where $i<j$ and the negative roots are those with $i>j$. The descent set of a permutation $g$ consists of the elements in $\Pi$ which $g$ maps to negative roots. For $L \subseteq \Pi$, let $X_L$ denote the set of permutations whose descent set is disjoint from $L$. It is well known \cite{H} that $|X_L|=\frac{n!}{|S_L|}$, where $|S_L|$ is the parabolic subgroup generated by the roots in $L$. Consequently if $p_L \geq 0$ satisfy $\sum_{L \subseteq \Pi} p_L =1$, the element $\sum_{L \subseteq \Pi} \frac{p_L |S_L|}{n!} X_L$ defines a probability measure on the symmetric group. 

	Given an element $\sum_{g \in S_n} c_g g$ of the group algebra of the symmetric group, by the inverse element we mean $\sum_{g \in S_n} c_g g^{-1}$. It is known that the RSK correspondence associates the same partition to $g$ and to $g^{-1}$, so when discussing the RSK correspondence one need not be concerned with whether we are considering an element in the group algebra or its inverse. The inverse of the element $\sum_{L \subseteq \pi} \frac{p_L |S_L|}{n!} X_L$ can be thought of as a shuffle. For instance if $p_{\Pi-\{\epsilon_1-\epsilon_2\}}=1$, this shuffle is simply the top to random shuffle. One reason these shuffles are important is a result of Solomon \cite{So} which states that $x_L x_K = \sum_{N \subseteq \Pi} a_{LKN} x_N$ for certain constants $a_{LKN}$. Thus one can at least in principle compute powers $(\sum_{L \subseteq \Pi} \frac{p_L |S_L|}{n!} X_L)^r$, which corresponds to understanding iterates of shuffles.  

	Now the main theorem of this section can be stated. Recall that the irreducible representations of the symmetric group $S_n$ are parameterized by partitions $\lambda$ of $n$.

\begin{theorem} \label{main2} For $L \subseteq \Pi$, let $J[L]$ denote the Markov chain associated to the pair $G=S_n$ and $H=S_L$, and let $J[\vec{p}]=\sum_L p_L J[L]$ denote the mixture of the Markov chains $J[L]$. Then $J[\vec{p}]_{1}^r(\lambda)$ (the chance that the mixed chain started at the trivial representation is at the representation parameterized by $\lambda$ after $r$ steps) is equal to the chance that an element of the symmetric group distributed as $(\sum_{L \subseteq \Pi} \frac{p_L |S_L|}{n!} X_L)^r$ has RSK shape $\lambda$. \end{theorem} 

\begin{proof} From Proposition \ref{spectrum}, the functions $\psi_C(\lambda)$ are a common orthonormal basis of eigenfunctions for the chains $J[L]$. Hence they are an orthonormal basis of eigenfunctions for the mixed chain $J[\vec{p}]$. This allows one to compute $J[\vec{p}]_{1}^r(\lambda)$ by the same method used in the proof of Theorem \ref{main1}, and one concludes that it is equal to 
\[ dim(\lambda) <\chi^{\lambda}, (\sum_L \frac{p_L |S_L|}{n!} Ind_{S_L}^{S_n}(1))^r>.\] 

	As explained in the preliminary remarks of Section 4 of \cite{BBHT}, the coefficients $a_{LKN}$ are related to tensor products of representations: \[ Ind_{S_L}^{S_n}(1) \times Ind_{S_K}^{S_n}(1) = \sum_{N \subseteq \Pi} a_{LKN} Ind_{S_N}^{S_n}(1).\] Letting $c_{N,r,\vec{p}}$ denote the coefficient of $X_N$ in \[ (\sum_{L \subseteq \Pi} \frac{p_L |S_L|}{n!} X_L)^r, \] it follows that $J[\vec{p}]_{1}^r(\lambda)$ is equal to \[ dim(\lambda) \sum_{N \subseteq \Pi} c_{N,r,\vec{p}} <\chi^{\lambda}, Ind_{S_N}^{S_n}(1)>.\] Letting $\mu$ denote the type of $N$ (that is $S_N$ is the direct product of symmetric groups whose sizes are the parts of the partition $\mu$), the multiplicity of $\lambda$ in $Ind_{S_N}^{S_n}(1)$ is by definition the Kostka-Foulkes number $K_{\lambda \mu}$ discussed in \cite{Sa}. Thus $J[\vec{p}]_{1}^r(\lambda)$ is equal to \[ dim(\lambda) \sum_{\mu} K_{\lambda,\mu} \sum_{N: type(N)=\mu} c_{N,r,\vec{p}} \] where the sum is over all partitions $\mu$ of $n$.

	Next it is necessary to show this is equal to the chance that an element of the symmetric group distributed as $(\sum_{L \subseteq \Pi} \frac{p_L |S_L|}{n!} X_L)^r$ has RSK shape $\lambda$. By the definition of $c_{N,r,\vec{p}}$, we know that \[ (\sum_{L \subseteq \Pi} \frac{p_L |S_L|}{n!} X_L)^r = \sum_{N \subseteq \Pi} c_{N,r,\vec{p}} X_N.\] So it suffices to show that the number of summands of the element $X_N$ (or equivalently the inverse of $X_N$) which the RSK correspondence maps to $\lambda$ is $dim(\lambda) K_{\lambda,type(N)}$. But writing $S_N=S_{a_1} \times S_{a_2} \cdots \times S_{a_r}$ the summands of the inverse of $x_N$ correspond (in an RSK shape preserving way) to words on the letters $\{1,\cdots,r\}$ in which the letter $l$ appears $a_l$ times. But such words with RSK shape $\lambda$ correspond to pairs $(P,Q)$ of Young tableau with $Q$ standard of shape $\lambda$ and $P$ semistandard of shape $\lambda$ and content $type(N)$. Since the number of these is $dim(\lambda) K_{\lambda,type(N)}$, the theorem is proved. \end{proof}

	Corollary \ref{totvar} is an important consequence of Theorem \ref{main2}. 

\begin{cor} \label{totvar} Let $tv(r,\vec{p})$ denote the total variation distance between the probability measure $(\sum_{L \subseteq \Pi} \frac{p_L |S_L|}{n!} X_L)^r$ on the symmetric group and the uniform distribution on the symmetric group. Let $\pi$ be the Plancherel measure of $S_n$. Then \[ \frac{1}{2} \sum_{\lambda \in Irr(S_n)} |dim(\lambda) <\chi^{\lambda}, (\sum_{L \subseteq \Pi} \frac{p_L |S_L|}{|n!|} Ind_{S_L}^{S_n}(1))^r> - \pi(\lambda)| \leq  tv(r,\vec{p}).\] \end{cor}

\begin{proof} From the proof of Theorem \ref{main1}, we know that 
\[ \frac{1}{2} \sum_{\lambda \in Irr(S_n)} |dim(\lambda) <\chi^{\lambda}, (\sum_{L \subseteq \Pi} \frac{p_L |S_L|}{|n!|} Ind_{S_L}^{S_n}(1))^r> - \pi_{\lambda}| \] is equal to the total variation distance between the measure $J[\vec{p}]_{1}^r$ and the Plancherel measure of the symmetric group. Theorem \ref{main2} gives that this is equal to the total variation distance between the RSK pushforward of the measure $(\sum_{L \subseteq \Pi} \frac{p_L |S_L|}{n!} X_L)^r$ and the Plancherel measure. Since the Plancherel measure is the RSK pushforward of the uniform distribution on the symmetric group, the corollary follows. \end{proof}

	The significance of Corollary \ref{totvar} is that it allows one to apply work on convergence rates of shuffles to the study of tensor products. We now give some examples which show that the bound of Corollary \ref{totvar} can be much sharper than that of Theorem \ref{main1}. 

{\bf Example 1: The defining representation}

	The first example is when $p_L=1$ for $L=\Pi-\{\epsilon_1-\epsilon_2\}$. Then $G=S_n$ and $H=S_{n-1}$. The representation theory problem in this case is the study of decompositions of the $rth$ tensor power of the defining (n-dimensional) representation, and the card shuffling problem is the $r$ fold iteration of the top to random shuffle.

	Consider the bound of Theorem \ref{main1}. Clearly $\beta=1-\frac{2}{n}$. It follows that \[ \sum_{\lambda \in Irr(S_n)} |\frac{dim(\lambda)}{n^r} <\chi^{\lambda}, (Ind_{H}^G(1))^r> - \pi(\lambda)| \leq \sqrt{n!}  (1-\frac{2}{n})^r.\] Using Stirling's approximation for $n!$ \cite{Fe}, one sees that for $r > \frac{nlog(n)+2c}{2 log(\frac{1}{\beta})}$, this is at most $(2 \pi)^{1/4} e^{-c}$. For $c$ fixed and large $n$, $\frac{nlog(n)+2c}{2 log(\frac{1}{\beta})}$ is roughly $\frac{n^2 log(n)}{4}$. 

	The bound from Corollary \ref{totvar} is much sharper. Indeed, it is known (\cite{AD}) that for $r=nlog(n)+cn$, the total variation distance between $r$ iterations of the top to random shuffle and the uniform distribution is at most $e^{-c}$, for $c \geq 0, n \geq 2$. 

	Next let us consider lower bounds for \[ \frac{1}{2} \sum_{\lambda \in Irr(S_n)} |\frac{dim(\lambda)}{n^r} <\chi^{\lambda}, (Ind_{H}^G(1))^r> - \pi(\lambda)| .\] By Theorem \ref{main2}, this is equal to the total variation distance between the RSK pushforward of $r$ iterations of the top to random shuffle and the Plancherel measure. A result of Chapter 5 of \cite{U} is that for large $n$ at least $\frac{1}{2} n log(n)$ iterations of the top to random shuffle are needed to randomize the length of the longest increasing subsequence (actually he states the result for the random to top shuffle, but this is the inverse of top to random). Since the longest increasing subsequence is a function of the RSK shape, it follows that \[ \frac{1}{2} \sum_{\lambda \in Irr(S_n)} |\frac{dim(\lambda)}{n^r} <\chi^{\lambda}, (Ind_{H}^G(1))^r> - \pi(\lambda)| \] requires $r$ at least $\frac{1}{2} n log(n)$ to be small. Thus the upper bound on $r$ in the previous paragraph is sharp to within a factor of two.

	The next two examples generalize example 1, but in different directions.

{\bf Example 2: $S_{n-k} \subset S_n$}

	This example is the case that $L=\Pi-\{\epsilon_1-\epsilon_2, \cdots, \epsilon_k-\epsilon_{k+1}\}$ where $k \leq n-1$. The representation theory problem is to study the decomposition of the $r$th tensor power of $Ind_{S_{n-k}}^{S_n}(1)$, and the relevant card shuffling is the top k to random shuffle, which proceeds by removing the top k cards from the deck and sequentially inserting them into random positions.

	First consider the bound of Theorem \ref{main1}. Using the
 fact that two elements in a symmetric group are conjugate if and only if
 they have the same structure, and that a conjugacy class with $n_i$
 cycles of length $i$ for all $i$ has size $\frac{n!}{\prod_i
 i^{n_i}n_i!}$, one finds that $\beta=\frac{(n-k)(n-k-1)}{n(n-1)}$. By
 the same argument as example 1, it follows that \[ \sum_{\lambda \in
 Irr(S_n)} |\frac{dim(\lambda)}{(n \cdots (n-k+1))^r} <\chi^{\lambda}, 
 (Ind_{H}^G(1))^r> - \pi(\lambda)| \leq (2 \pi)^{1/4} e^{-c}\]
 when $r > \frac{nlog(n)+2c}{2 log(\frac{1}{\beta})}$. For fixed $c,k$
 and large $n$, $\frac{nlog(n)+2c}{2 log(\frac{1}{\beta})}$ is roughly
 $\frac{n^2 log(n)}{4k}$.

	The convergence rate of the card shuffling problem was studied in \cite{DFP}, where it was shown that for $k$ fixed and large $n$, the total variation distance is at most $e^{-c}$ for $r=\frac{n}{k}(log(n)+c)$. Thus the bound from Corollary \ref{totvar} is much sharper. The argument for the lower bound also generalizes, showing that $r$ must be at least $\frac{1}{2k} n log(n)$ for \[ \frac{1}{2} \sum_{\lambda \in Irr(S_n)} |\frac{dim(\lambda)}{(n(n-1) \cdots (n-k+1))^r} <\chi^{\lambda}, (Ind_{H}^G(1))^r> - \pi(\lambda)| \] to be small.

{\bf Example 3: Action on k-sets}

	The next example is the case that $p_L=1$ where
	$L=\Pi-\{\epsilon_k-\epsilon_{k+1}\}$, and $1 \leq k \leq
	n/2$. Then $G=S_n$ and $H=S_k \times S_{n-k}$. The
	representation theory problem in this case is the study of
	decompositions of the $rth$ tensor power of the permutation
	representation on $k$-sets, and the card shuffling problem is
	the $r$ fold iteration of the shuffle which proceeds by
	cutting off exactly $k$ cards, and then riffling them with the
	other $n-k$ cards (i.e. choosing a random interleaving).

	First consider the bound of Theorem \ref{main1}. The value of $\beta$ is calculated in \cite{GM} for $n \geq 5$ and shown to occur for the conjugacy class of transpositions, where it is $\frac{{n-2 \choose k} + {n-2 \choose k-2}}{{n \choose k}}$. For $k$ fixed and large $n$, $log(\frac{1}{\beta})$ is roughly $\frac{2k}{n}$, so that $r$ slightly more than $\frac{n^2 log(n)}{4k}$ will make \[ \sum_{\lambda \in Irr(S_n)} |\frac{dim(\lambda)}{{n \choose k}^r} <\chi^{\lambda}, (Ind_{H}^G(1))^r> - \pi(\lambda)| \] small.

	Now consider the bound from Corollary \ref{totvar}. To apply
it we require an upper bound on the total variation distance between
the uniform distribution and $r$ iterations of the shuffle which cuts
off exactly $k$ cards and riffles them with the rest of the deck. This
shuffle is a special case of the Bidigare-Hanlon-Rockmore walks on
chambers of hyperplane arrangements, and a convenient upper bound for total
variation distance is in \cite{BD} (this bound is somewhat weaker than
the bound in \cite{BHR} but is easier to apply). In the case at hand
one can check that the total variation bound becomes ${n \choose 2}
\left( \frac{{n-2 \choose k-2} + {n-2 \choose k}}{{n \choose k}}
\right)^r$, which is better than the bound $\frac{\sqrt{n!}}{2}
\left( \frac{{n-2 \choose k} + {n-2 \choose k-2}}{{n \choose k}}
\right)^r$ from Theorem \ref{main1}. One concludes that $r$ slightly
more than $\frac{n log(n)}{k}$ makes \[ \frac{1}{2} \sum_{\lambda \in
Irr(S_n)} |\frac{dim(\lambda)}{{n \choose k}^r} <\chi^{\lambda}, 
(Ind_{H}^G(1))^r> - \pi(\lambda)| \] small. Moreover, the argument for the lower bound in the other examples generalizes, showing that $r$ must be at least $\frac{n log(n)}{2k}$.

	We remark that the fact that non-reversible Markov chains such as top to random are related to the reversible Markov chain $J$ by means of Theorem \ref{main2} is quite mysterious. As a further result in this direction, we show that the Markov chains $J[\vec{p}]$ and $\sum_{L \subseteq \Pi} \frac{p_L |S_L|}{n!} X_L$ are isospectral. See \cite{F4} for some other connections between the top to random shuffle and nonreversible Markov chains.

\begin{prop} The Markov chain $J[\vec{p}]$ and the element $\sum_{L \subseteq \Pi} \frac{p_L |S_L|}{n!} X_L$ have the same set of eigenvalues. \end{prop}

\begin{proof} Since the chains $J[L]$ have a common basis of eigenvectors, the eigenvalues of $J[\vec{p}]$ are linear functions in the $p$'s. Similarly \cite{BHR} finds a formula for the eigenvalues of the element $\sum_{L \subseteq \Pi} \frac{p_L |S_L|}{n!} X_L$ and shows that they are linear in the $p$'s. Hence it is enough to prove the result when $p_L=1$ for some $L$. 

	From Corollary 2.2 of \cite{BHR}, the eigenvalues of the
element $\frac{|S_L|}{n!} X_L$ are indexed by permutations $g \in
S_n$. Let $\mu$ be such that the orbits of $S_L$ on $\{1,\cdots,n\}$
are $\{1 \cdots \mu_1\}, \{\mu_1+1 \cdots \mu_1+\mu_2 \}$,
etc.; hence $\mu$ is a composition of $n$. A block ordered partition of the set $\{1,\cdots,n\}$ is by definition a set partition with an ordering on the blocks of the partition. We say that a block ordered partition has type $\mu$ if the first block has size $\mu_1$, the second block has size $\mu_2$ and so on. The result of \cite{BHR} is that the eigenvalue corresponding to $g$ is the proportion of block ordered partitions of type $\mu$ which are fixed by $g$ in the sense that each block is sent to itself. This is equivalent to requiring that each block is a union of cycles of $g$. Letting $n_i$ denote the number of $i$-cycles of $g$, it follows that this proportion is \[ \frac{\mu_1! \mu_2! \cdots}{n!} \sum_{\sum_k a_i^{(k)}= n_i \atop \sum_{i} ia_i^{(k)}=\mu_k \ all \ k} \prod_{i \geq 1} {n_i \choose a_i^{(1)},a_i^{(2)},\cdots}.\]  

	On the other hand, by Proposition \ref{spectrum}, we know that the eigenvalues of $J$ are parameterized by conjugacy classes $C$ of $S_n$. Let $n_i$ denote the number of cycles of length $i$ for elements in the class $C$. Using the fact that $|C|=\frac{n!}{\prod_i i^{n_i}n_i!}$, it follows that \[ \frac{|C \cap S_L|}{|C|}  = \frac{\prod_i i^{n_i}n_i!}{n!} \sum_{\sum_k a_i^{(k)}= n_i \atop \sum_{i} ia_i^{(k)}=\mu_k \ all \ k} \prod_k \frac{\mu_k!}{ \prod_i i^{a_i^{(k)}} a_i^{(k)}!}. \] This is equal to the expression of the previous paragraph, so the proof is complete. \end{proof}

\end{document}